\input amstex
\documentstyle{amsppt}
\magnification=\magstep 1
\voffset-0.3cm
\TagsOnLeft
\loadbold
\topmatter
\title  Linking and coincidence invariants   \endtitle
\author Ulrich Koschorke \\ { }\endauthor
\leftheadtext{Ulrich Koschorke}
\address Universit\"at Siegen,
Emmy Noether Campus, Walter-Flex-Str. 3,
D-57068 Siegen, Germany
\endaddress
\email koschorke\@mathematik.uni-siegen.de \endemail
\abstract
Given a link map \ $f$ \ into a manifold of the form \
$Q = N \times \Bbb R$, \ when can it be deformed to an \lq\lq
unlinked\rq\rq \ position (in some sense, e.g. where its components map to
disjoint \ $\Bbb R$-levels)? Using the language of normal bordism
theory as well as the path space approach of Hatcher and Quinn we
define obstructions \ $\widetilde\omega_\varepsilon (f), \
\varepsilon = +$ \ or \ $\varepsilon = -$, \ which often answer
this question completely and which, in addition, turn out to
distinguish a great number of different link homotopy classes.
In certain cases they even allow a complete link homotopy classification.

Our development parallels recent advances in Nielsen coincidence
theory and leads also to the notion of Nielsen numbers of link
maps.

In the special case when \ $N$ \ is a product of spheres sample
calculations are carried out. They involve the homotopy theory of
spheres and, in particular, James--Hopf--invariants.
\endabstract

 \subjclass Primary 55P35, 55S35, 57Q45,
57R90. Secondary 55M20, 55Q25, 55Q45
\endsubjclass

\endtopmatter

\input BoxedEPS.tex
\SetRokickiEPSFSpecial
\HideDisplacementBoxes

\input xy
\input xymatrix
\input xyarrow
\input xycurve

\define\incl{\operatorname{incl}}
\define\sign{\operatorname{sign}}

\define\scirc{{\ssize{\circ}}}

\document

\specialhead 1.\ \ Introduction
\endspecialhead

Throughout this paper let \ $M^{m_1}_1, M^{m_2}_2, N^n$ \ and \
$Q$ \ denote smooth manifolds (of the indicated dimensions)
without boundary, where \ $M_1, M_2$ \ are compact and \ $N$ \ is
connected, and let
$$
f \ \ = \ \ f_1 \ \amalg \ f_2 \ \ : \ \ M_1 \ \amalg \ M_2 \ \
\longrightarrow \ \ Q
$$
be a \ {\it link map} \ (i.e.\ the continuous maps \ $f_1$ \ and \
$f_2$ \ have disjoint images).

Two such link maps are called \ {\it link homotopic} \ (compare
Milnor \cite{M}) if they can be deformed continuously into one
another through a family of link maps.

\subhead Question
\endsubhead \
When can \ $f$ \ be unlinked? More precisely: \ when is \ $f$ \
link homotopic to a link map which is trivial in some sense?
\medskip

In classical link theory two approaches to such problems have
played a central role: consider either appropriate intersections
or (over)crossings.

If there is some canonical notion of what a  \ \lq\lq
trivial\rq\rq \ or \ \lq\lq faraway\rq\rq\ position of \ $f_1$ \
should be, and if a homotopy \ $F_1$ moves \ $f_1$ \ to such a
position, measure the \ {\it intersection} \ locus \ $C$ \ of \
$F_1$ \ with \ $f_2$ \ in some way. (E.g.\ if the domain of \
$f_1$ \ is a sphere and \ $F_1$ \ is a nulhomotopy this approach
leads to the standard procedure of intersecting \ $f_2$ \ with a
singular ball spanned by \ $f_1$.)

If \ $Q$ \ has the special product form \ $Q = N \times \Bbb R$ \
there are natural choices \ $F_+$ \ and \ $F_-$ \ for such a
homotopy: \ we can move \ $f_1$ in the positive (or negative) \
$\Bbb R$-direction until the whole image of \ $f_1$ \ lies above
(or below) the image of \ $f_2$ \ w.r.\ to the \ $\Bbb R$-levels.
The intersection of such a homotopy with \ $f_2$ corresponds to
the \ {\it overcrossing} \ (or \ {\it undercrossing}, resp.) locus
\ $C_{\pm}$ \ of the projections \ $f'_1$ \ and \ $f'_2$ \ to \
$N$.

Whether we base our approach on intersections or
over/undercrossings, the resulting unlinking obstruction will be
all the more powerful if we register all relevant geometric data
concerning the locus \ $C$ \ or \ $C_{\pm}$ \ as carefully as
possible. One rather obvious strategy is to use the language of
normal bordism theory. It keeps track of the relations between the
stable tangent or normal bundles of the intersection or crossing
locus on one hand, and of the manifolds \ $M_1, M_2,$ \ and \ $N$
\ on the other hand. In very special situations this amounts just
to framed bordism (involving stably parallelized manifolds), but
in general normal bordism is much more widely applicable and
flexible, and a much stronger tool than e.g.\ oriented bordism (if
it applies) or (co~)-homology with twisted coefficients.

An additional refinement was inspired by the fundamental work of
Hatcher and Quinn \cite{HQ}. It is easily overlooked that every
coincidence point \ $x$ \ comes naturally with a path, namely the
constant path at the common value \ $f'_1 (x) = f'_2 (x)$. But
this datum carries very valuable information. Keeping track of it
and accommodating our obstruction accordingly in a normal bordism
group of a suitable path space will in certain situations supply
the necessary data needed to construct a homotopy which deforms
maps away from one another or which unlinks link maps.

In sections 2 and 3 of this paper we define and study unlinking
obstructions
$$
\widetilde\omega_+ (f) \ , \  \widetilde\omega_- (f) \ \ \in \ \
\Omega_{m_1 + m_2 - n} (E(f); \ \widetilde\varphi)
$$
which often lead to a complete answer of our original question
(see theorem 2.13 below). They lie in a normal bordism group of an
appropriate pathspace \ $E (f)$. In many interesting cases this
space has an extremely rich topology. E.g.\ already the set \
$\pi_0 (E (f))$ \ of pathcomponents may be  huge (there is a
natural bijection onto a certain well-studied quotient of \ $\pi_1
(N)$, the so called Reidemeister set). The resulting decomposition
of our invariants allows us to define \ {\it Nielsen numbers} \
$N_+ (f)$ \ and \ $N_- (f)$ \ of a link map: \ just count the
(finitely many) \ {\it essential} \ pathcomponents of \ $E (f)$, \
i.e.\ those where the corresponding components of \
$\widetilde\omega_{\pm} (f)$ \ are nontrivial. This procedure
replaces the often unwieldy \
$\widetilde\omega_{\varepsilon}$-obstruction (which e.g.\ lies in
a group varying with \ $f \ $) by the numerical link homotopy
invariant \ $N_{\varepsilon}$ \ which vanishes precisely if \
$\widetilde\omega_{\varepsilon}$ \ does, \ $\varepsilon = +$ or
$-$.

A similar point of view was recently introduced into the study of
fixed point and coincidence phenomena and lead to a coherent
Nielsen coincidence theory involving manifolds with arbitrary
orientation behaviour and dimension combinations (cf. \cite{K 5}).

Our approach is also closely related to recent work of A.\ Pilz
(cf.\ \cite{P}). His invariant \ $\alpha_w (f)$ \ registers the
decomposition of the bordism class of the intersection \ $F_+
\pitchfork f_2$ \ (in framed, oriented or unoriented bordism \;
$\Omega^{fr}_*, \ \Omega_*$ \ or \ $\frak N$ \ -- as the situation
may permit) into components indexed by the Reidemeister set. This
is strong enough to yield a link homotopy classification result
when \ $m_1 + m_2 = n$ \ (cf.\ \cite{P}, theorem 5.4).

If \ $M_1$ \ and \ $M_2$ \ are sufficiently highly connected then
natural isomorphisms (exhibited by Hatcher and Quinn) allow us to
interpret \ $\widetilde\omega_{\pm}$ \ itself as a link homotopy
invariant. It takes values in the \ $(m_1 + m_2 - n)$th framed
bordism group of the loop space \ $\Lambda N$ \ of \ $N$.

In particular, this applies when \ $M_1$ \ and \ $M_2$ \ are
spheres of dimensions \ $\le n -2$. In this case there is also a
well-defined addition of link maps \ and \
$\widetilde\omega_{\pm}$ \ turns out to be compatible with this
and other natural operations (cf.\ 5.2 - 5.4 below). Moreover a
simple construction (using \lq\lq meridians\rq\rq) supplies many
examples of link maps with interesting \
$\widetilde\omega_\pm$-values (cf.\ 5.7 - 5.11). Thus our
invariants -- originally conceived as unlinking obstructions --
turn out to distinguish also a great number of different link
homotopy classes. In some situations they even lead to a complete
classification.

\proclaim{Theorem 1.1} \ Assume \ $1 \le m_1 + 1, m_2 \le 2n - m_1
- m_2 - 2$ \ and that \ $N$ \ is stably parallelized.

Then two base point preserving link maps
$$
f\ , \ \widehat f \ \ \ : \ \ \ S^{m_1} \amalg S^{m_2} \
\longrightarrow \ N \times \Bbb R
$$
are link homotopic $($in the base point preserving sense$)$ if and
only if \ $[f_i] = [\widehat f_i] \in \pi_{m_i} (N \times \Bbb R),
\ i = 1, 2$, \ and if in addition
$$
\widetilde\omega_+ (f) \ = \ \widetilde\omega_+ (\widehat f) \ \
\in \ \ \Omega^{fr}_{m_1 + m_2 - n} (\Lambda N)
$$
$($or, equivalently, \ $\widetilde\omega_- (f) =
\widetilde\omega_- (\widehat f))$.
\endproclaim

This is proved in section 5; the relation between base point
preserving and base point free link homotopy theory is indicated
in remarks 3.9 and 4.2.

For an illustration we consider the case when \ $N$ \ is a product
of spheres. In section 4 we reduce the calculation of the framed
bordism groups of the loop space \ $\Lambda N$ \  (which
accommodate \ $\widetilde\omega_\pm (f)$) \ via
James-Hopf-invariants to standard methods of the stable homotopy
theory of spheres. This can be used in many concrete settings such
as

\example{Example 1.2: \ $\bold N \bold= \bold S^1 \bold\times
\bold S^2 \bold\times \bold S^8$} \ The (base point {\it free})
link homotopy class of a link map
\vskip1mm
$$
f = f_1 \amalg f_2 \ : \ S^6 \amalg S^7 \ \longrightarrow S^1
\times S^2 \times S^8 \times  \Bbb R
$$
\vskip1mm
is completely determined by the homotopy classes
$$
[f_1] \ \in \ [S^6, N \times \Bbb R] \ \cong \ \Bbb Z_{12} \ \ , \
\ [f_2]  \ \in \  [S^7, N \times \Bbb R] \ \cong \ \Bbb Z_2
$$
of the component maps and by the unlinking obstruction \
$[\widetilde\omega_+ (f)]$ \ which lies in
$$
\Omega^{fr}_2 (\Lambda N)/\!\!\sim \ \ \cong \ \
\Bigl(\bigoplus^\infty_{- \infty} (\Omega^{fr}_2 \oplus
\Omega^{fr}_1 \oplus \Omega^{fr}_0)\Bigr)/\!\!\sim \ \ \cong \ \
(\Bbb Z_2 \oplus \Bbb Z_2 \oplus \Bbb Z)[X^{\pm 1}]/\!\!\sim .
$$
Here two (formal) Laurent polynomials (with coefficients in the
group \ $\Bbb Z_2 \oplus \Bbb Z_2 \oplus \Bbb Z$) are equivalent
if they differ by the factor $X^j$\ for some integer $j$\.
\newline This follows from theorems 1.1, 4.1, remarks 3.9, 4.2 and the
tables of Toda \cite{T}.
\endexample

\example{Example 1.3: $\bold N \bold= (\bold S^1)^3 \bold\times
\bold S^2$} \ Both for \ $\varepsilon = +$ \ and \ $\varepsilon =
-$ \ every element of
$$
\Omega_1^{fr} (\Lambda ((S^1)^3 \times S^2)) \ \ \ \cong \ \ \
\bigoplus_{j \in \Bbb Z^3} (\Bbb Z_2 \oplus \Bbb Z)
$$
occurs as the \ $\widetilde\omega_\varepsilon$--value of a (base
point preserving) link map
$$
f \ : \ S^3 \amalg S^3 \ \longrightarrow \ (S^1)^3 \times S^2
\times \Bbb R \ .
$$
This and further examples will be discussed at the end of section
5.
\endexample

\vskip7mm

\specialhead  2.\ \ The unlinking obstructions \
$\widetilde\omega_{\pm}$ \ and \ $\omega_{\pm}$
\endspecialhead

In this section we adapt the coincidence invariants \
$\widetilde\omega$ \ and \ $\omega$ \ constructed in \ \cite{K 5}
to the setting of link maps. We will obtain the obstructions \
$\widetilde\omega_{\varepsilon}$ \ and \ $\omega_{\varepsilon}$ \
where  \ $\varepsilon$ \ stands for the symbols \ $+$ \ or \ $-$ \
(and, in formulas, for the factors \ $+ 1$ \ and \ $- 1$).

Throughout the remainder of this paper (unless mentioned
otherwise)
$$
f \ \ = \ \ f_1 \ \amalg \ f_2 \ \ : \ \ M_1 \ \amalg \ M_2 \ \
\longrightarrow \ \ N \times \Bbb R \tag 2.1
$$
will be a link map into a manifold of the indicated product form,
and \ $f_i = (f'_i, f''_i)$ \ is the corresponding decomposition
via the projections to \ $N$ \ and \ $\Bbb R$, resp., \ $i = 1,
2$. Consider also the product manifold \ $M := M_1 \times M_2$ \
with the dimension \ $m := m_1 + m_2$ \ and the projections \ $p_i
: M \to M_i, \ i = 1,2$.

Our discussion will center around the space \footnote"*"{In
\cite{HQ} (and \cite{K 5}, resp.) this space is denoted by $E
(f'_1, f'_2)$ (and by \newline $E (f'_1 \scirc p_1, f'_2 \scirc
p_2$), resp.)}
$$
E (f) \ := \ \{ ((x_1, x_2); \theta) \ \in \ M \times P (N) \ \ |
\ \  \theta (0) = f'_1 (x_1), \ \theta (1) = f'_2 (x_2)\} \tag 2.2
$$
where \ $P (N)$ \ is the set of all continuous paths \ $\theta :
[0, 1] \to N$, \ endowed with the compact -- open topology. Let \
$pr : E (f) \to M$ \ denote the obvious fiber projection.

If the map
$$
(f'_1 \scirc p_1, f'_2 \scirc p_2) \ = \ f'_1 \times f'_2 \ : \ M
= M_1 \times M_2 \ @>{ \ \ \ }>> \ N \times N
$$
is smooth and transverse to the diagonal
$$
\Delta \ := \ \{ (y, y) \ \in \ N \times N \ | \ y \in N \} \tag
2.3
$$
then the \ $\varepsilon$--{\it coincidence locus} \
$$
C_{\varepsilon} (f) \ := \ \{ (x_1, x_2) \in M_1 \times M_2 \ | \
f'_1 (x_1) = f'_2 (x_2), \ \varepsilon f^{''}_1 (x_1) <
\varepsilon f^{''}_2 (x_2)\} \tag 2.4
$$
(where \lq\lq $f_2$ \ {\it over}crosses \ $f_1$\rq\rq\  if \
$\varepsilon = +$ \ and \  $f_2$ \ {\it under}crosses \ $f_1$ \ if
\ $\varepsilon = -$) \ is a closed smooth \ $(m - n)$-dimensional
manifold, equipped with the map
$$
\widetilde g_{\varepsilon} \ : \ C_{\varepsilon} (f)  @>{ \ \ \
}>> \ E (f) \tag 2.5
$$
(which sends \ $(x_1, x_2) \in C_\varepsilon (f)$ \ to \ $((x_1,
x_2)$, constant path at \ $f'_1 (x_1) = f'_2 (x_2))$ \ and with a
stable vector bundle isomorphism
$$
\overline g_{\varepsilon} \ : \ TC_\varepsilon (f) \oplus
\widetilde g^*_{\varepsilon} (pr^* ((f'_1 \scirc p_1)^* (TN))) \ \
\cong \ \ \widetilde g^*_{\varepsilon} (pr^* (TM)) \tag 2.6
$$
(since the normal bundle \ \ $\nu (C_{\varepsilon} (f), M)$ \ \ is
canonically isomorphic to \ \ $(f'_1 \scirc p_1)^* (TN))$.

If \ $f$ \ is an arbitrary link map, apply this construction to an
approximation of \ $f'_1 \times f'_2$ \ which is smooth and
transverse to \ $\Delta$.

In any case the resulting triples \ $(C_+, \widetilde g_+,
\overline g_+)$ \ and \ $(C_-, \widetilde g_-, \overline g_-)$ \
determine well-defined normal bordism classes
$$
\widetilde\omega_+ (f) \ , \ \widetilde\omega_- (f) \ \in \
\Omega_{m - n} (E (f); \widetilde\varphi) \tag 2.7
$$
as in \cite{K 5}, \S\ 4; here
$$
\aligned
\varphi \ := \ &(f'_1 \scirc p_1)^* (TN) \ - \ TM  \ \qquad \qquad \
\text{and} \\
\widetilde\varphi \ := \ &pr^* (\varphi)
\endaligned
\tag 2.8
$$
are the relevant (virtual) coefficient bundles over \ $M$ \ and \
$E(f)$, resp.

Clearly we have
$$
\widetilde\omega (f'_1 \scirc p_1, f'_2 \scirc p_2) \ = \
\widetilde\omega_+ (f) + \widetilde\omega_- (f) \tag 2.9
$$
where the left hand term is the full coincidence invariant
discussed in \cite{K 5}. Indeed, all we have done here is to
decompose the coincidence locus \ $(f'_1 \times f'_2)^{- 1}
(\Delta)$ \ disjointly into its overcrossing and undercrossing
parts.

If we forget the path space aspect of our data and keep track only
of the over- or undercrossing manifolds, together with the way
they sit in \ $M$ \ and with their \lq\lq twisted framings\rq\rq \
$\overline g_{\varepsilon}$, \ we obtain the weaker invariants
$$
\omega_{\varepsilon} (f) \ := \ [\, C_{\varepsilon}\,, \
\text{inclusion}\,, \ \overline g_{\varepsilon}\, ] \ \ = \ \ pr_*
(\widetilde\omega_{\varepsilon} (f)) \ \in \ \Omega_{m +n} (M;
\varphi) \ . \tag 2.10
$$
\medskip

\definition{Definition 2.11} \
The link map \ $f$ \ is \ $\varepsilon$--{\it unlinked} \ if
$$
\varepsilon (f^{''}_1 (x_1) - f^{''}_2 (x_2)) \ \ > 0 \ \ \
\text{for all} \ (x_1, x_2) \in M_1 \times M_2
$$
(compare 2.1) or, equivalently, if the image \ $f_1 (M_1)$ \ in \
$N \times \Bbb R$ \ lies strictly above (or below, resp.) \ $f_2
(M_2)$ \ w.r.\ to the \ $\Bbb R$--coordinate when \ $\varepsilon =
+$ \ (or $\varepsilon = -$, resp.).

$f$ \ is called \ {\it $\varepsilon$-unlinkable} \ if \ $f$ \ is
link homotopic to an \ $\varepsilon$-unlinked link map.
\enddefinition
\medskip

\proclaim{Proposition 2.12} \ If \ $f$ \ is \
$\varepsilon$--unlinkable, then \ $\widetilde\omega_{\varepsilon}
(f) = 0$ \ and therefore also \ $\omega_{\varepsilon} (f) = 0$.
\endproclaim

\demo{Proof} \ If \ $f$  is  \ $\varepsilon$--unlinked, then \
$C_{\varepsilon} (f) = \phi$ \ and hence \
$\widetilde\omega_{\varepsilon} (f) = 0$. Moreover, recall from
\cite{K 5}, \S\ 3--4, that homotopies of \ $f_i, \ i = 1, 2$, \
induce isomorphisms of normal bordism groups which are compatible
with the full coincidence invariants \ $\widetilde\omega$ \ and --
if we are dealing with \ {\it link} \ homotopies -- preserve even
the decomposition \ $\widetilde\omega = \widetilde\omega_+ +
\widetilde\omega_-$. \hfill $\blacksquare$
\enddemo
\medskip

The methods of Hatcher and Quinn yield the following converse
result.

\proclaim{Theorem 2.13} \ Assume that \  $m_1 + 1, m_2 \le 2n -
m_1 - m_2 - 2$ \ or \ $m_1, \ m_2 + 1 \le 2n - m_1 - m_2 - 2$.

If after a suitable link homotopy \ $f_1$ \ and \ $f_2$ \ project
to smooth immersions into \ $N$ \ $($this holds in particular if \
$M_1, \ M_2$ \ and \ $N$ \ are stably parallelisable$)$, then we
have for \ $\varepsilon = +$ \ and \ $\varepsilon = -$:

$f$ is \ $\varepsilon$-unlinkable precisely when \
$\widetilde\omega_\varepsilon (f) = 0$.
\endproclaim

\demo{Proof} \ Our claim is valid for \ $f$ \ and \ $\varepsilon$
\ if and only if it holds for \ $f_2 \amalg f_1$ \ and \ $-
\varepsilon$ (compare the discussion of (4.5) in \cite {K 5}).
After possibly interchanging \ $f_1$ \ and \ $f_2$ \ we may
therefore assume the first of the above mentioned inequalities.

Apply the generalized Whitney trick construction of the proof of
theorem 2.2 in \cite{HQ} to the immersions \ $F_1 = F_\varepsilon$
\ (compare 2.14 below, or the beginning of our introduction) and \
$f_2$, as well as to a nulbordism of the over/undercrossing data
of \ $C_\varepsilon \approx F_1 \pitchfork f_2$. \ The resulting
deformation will move \ $F_1, f_2$ \ to maps \ $\widehat F_1,
\widehat f_2$ \ with disjoint images. Since the key steps of the
construction are based on approximations we can make sure that it
does not interfere with \ $F_1 \ | \ M_1 \times \{ 0, 1\}$ \ and
faraway \ $\Bbb R$-levels. Thus \ $f$ \ is link homotopic to \
$f_1 \amalg \widehat f_2$ \ and, via \ $\widehat F_1$, \ to an \
$\varepsilon$-unlinked link map. \hfill $\blacksquare$
\enddemo

\remark{Remark $2.14$} \ The previous proof is based on a
generally valid alternate interpretation of our invariants. Let
$$
F_1 \ = \ (f'_1, F^{''}_1) \ : \ M_1 \times [0, 1] \ @>{ \ \ \ }>>
\ N \times \Bbb R
$$
be a homotopy which deforms \ $F_1 (-, 0) = f_1$ \ monotonously in
the positive $($or negative$)$ \ $\Bbb R$--direction until \ $F_1
(-, 1) \amalg f_2$ \ is \ $\varepsilon-$unlinked, \ $\varepsilon =
+ $ \ (or \ $\varepsilon = -$, resp.). Then \ $F_1$ \ and \ $f_2$
\ define a pair of maps from \ $M_1 \times M_2 \times [0, 1]$ \ to
\ $N \times \Bbb R$ \ whose coincidence manifold is essentially \
$C_{\varepsilon} (f)$, \ with compatible coincidence data. Thus
the resulting normal bordism class corresponds to \ $\varepsilon
\widetilde\omega_{\varepsilon} (f)$ \ via the isomorphism induced
by the projections \ $M \times (0, 1) \sim M$ \ and \ $N \times
\Bbb R \sim N$.

This alternate (\lq\lq intersection\rq\rq) approach allows us
sometimes to extend our invariants to link maps into more general
target manifolds \ $Q$ \ $($e.g.\ if \ $M_1 = S^{m_1}$ \ and \
$\pi_{m_1} (Q) = \pi_{m_1 + 1} (Q) = 0\,)$  \hfill $\blacksquare$
\endremark
\medskip

Now consider the special case that \ $N$ \ has the form \ $N = N'
\times \Bbb R$. Then \ $\widetilde\omega (f'_1 \scirc p_1, f'_2
\scirc p_2) = 0$ \ since the extra \ $\Bbb R$--direction allows us
to move \ $f'_1 \scirc p_1$ \ and \ $f'_2 \scirc p_2$ \ apart.
Thus \ $\widetilde\omega_+ (f) = - \widetilde\omega_- (f)$ \ (cf.
2.9). In fact, we can use the projection along \ {\it any} \ ray \
$\Bbb R_+ \cdot v$ \ in \ $\Bbb R \times \Bbb R$ \ (where \ $v \ne
0$) \ to study the over/undercrossing behaviour of link maps into
\ $N' \times \Bbb R \times \Bbb R$.

\example{Example 2.15: \ Classical link maps} \ Here \ $M_i =
S^{m_i}, \ N = \Bbb R^n, \ n > 0$, \ and \ $f : S^{m_1} \amalg
S^{m_2} \to \Bbb R^{n + 1}$. Because of the linear structure on \
$\Bbb R^{n + 1}$ the fiber projection \ $pr : E (f) \to S^{m_1}
\times S^{m_2}$ \ is a homotopy equivalence so that \
$\widetilde\omega_+ (f) = - \widetilde\omega_- (f)$ \ is precisely
as strong as the invariant (cf. 2.10)
$$
\omega_+ (f) \ \in \ \Omega^{fr}_{m - n} (S^{m_1} \times S^{m_2})
\ \cong \ \Omega^{fr}_{m - n} \oplus \Omega^{fr}_{m_1 - n} \oplus
\Omega^{fr}_{m_2 - n}
$$
In this decomposition we use the isomorphism \ $\Omega^{fr}_*
(S^{m_i} \times X) \cong \Omega^{fr}_* (X) \oplus \Omega^{fr}_{* -
m_i} (X)$ \ defined by the projection to  $X$  and by transverse
intersection with  $\{*_i\} \times X;$ \ \ $\omega_+ (f)$
corresponds to the triple
$$
(\alpha (f), \ \ \omega_+ (f | (S^{m_1} \amalg \{ *_2\})), \ \
\omega_+ (f | (\{*_1\} \amalg S^{m_2})))
$$
consisting of the \lq\lq generalized linking number\rq\rq\ $\alpha
(f)$ \ (cf.\ \cite{K 2}) and of the overcrossing invariants of one
sphere with just the base point of the other sphere.

In the dimension range of theorem 2.13 the second and third
components of this triple vanish; thus \ $f$ \ is \
$\varepsilon$-unlinkable or, equivalently, link nulhomotopic
precisely if \ $\alpha (f) = 0$. \ Actually, N.\ Habegger and U.\
Kaiser \cite{HK} have shown that the \ $\alpha$-invariant
classifies \ $f$ \ completely up to link homotopy in the more
general range \ $2 (m_1 + m_2) \le 3n - 2, \ \ \ m_1, m_2 < n$.
(Compare also \cite{S}). \hfill $\blacksquare$
\endexample

In contrast to this example we will see below that \
$\widetilde\omega_{\varepsilon} (f)$ \ is often considerably
stronger than \ $\omega_{\varepsilon} (f)$.

\vskip7mm

\specialhead 3.\ \ Nielsen numbers of link maps and other link
homotopy invariants
\endspecialhead

In order to get a better understanding of our invariants and of
the groups in which they lie we need to recall a few facts about
the Hurewicz fibration \ $pr : E (f) \to M$ \ (compare \cite{K 5},
\S\ 2).

Pick points \ $x_0 = (x_{10}, x_{20}) \in M = M_1 \times M_2, \ \
y_0 \in N$ \ and paths \ $\sigma_i$ \ in \ $N$ \ joining \ $y_0$ \
to \ $f'_i (x_{i0}), \ i = 1, 2$. This choice determines a
homotopy equivalence between the space \ $\Lambda (N, y_0)$ \ of
loops in \ $N$ \ starting and ending in \ $y_0$, \ and the fiber \
$pr^{- 1} (x_0)$ \ which consists of all paths in \ $N$ \ from \
$f'_1 (x_{10})$ \ to \ $f'_2 (x_{20})$. We compose with the fiber
inclusion to obtain
$$
\incl \  : \ \Lambda (N, y_0) \ \longrightarrow \ E (f) \ . \tag
3.1
$$
\medskip

\proclaim{Proposition 3.2 (cf. \cite{K 5}, 2.1)} \ The induced map
$$
\incl_* \ \  : \  \ \pi_1 (N, y_0) \ = \ \pi_0 (\Lambda (N, y_0))
\  @>{ \ \ \ \ }>> \ \pi_0 (E (f))
$$
yields a bijection from the so called  \ {\rm Reidemeister set}
\vskip2mm
$$
R \ \ \ := \ \ \ \xy
   \POS(18.7,-3)*+{\sigma_{1*} f'_{1*}(\pi_1(M_1,x_{10}))},
      \POS(38.2,-6) \ar@{-} (29.2,6),
    \POS(41,3)*+{\pi_1(N,y_0)},
      \POS(43.8,-6) \ar@{-} (52.8,6),
    \POS(63.2,-3)*+{\sigma_{2*}f'_{2*}(\pi_1(M_2,x_{20}))},
\endxy
$$
\vskip2mm
\flushpar
 onto the set of pathcomponents of \ $E (f)$ \ (Here \
$\sigma_{i*}: \pi_1 (N, f'_i (x_{i0}))\to \pi_1 (N, y_0)$ \
denotes the obvious isomorphism induced by \ $\sigma_i, \ i = 1,
2)$.
\endproclaim

Thus the rich geometry of \ $E (f)$ \ manifests itself already in
a possibly very large number of pathcomponents. However, since the
coincidence manifold \ $C_{\varepsilon}$ \ is compact, only
finitely many pathcomponents \ $A \in \pi_0 (E (f))$ \ are \ {\it
essential}, \ i.e.\ the corresponding direct summand
$$
\widetilde\omega_{\varepsilon, A} (f) \ = \ [ \ C_{\varepsilon, A}
(f)  :=  \widetilde g^{- 1} (A), \ \ \widetilde g |
C_{\varepsilon, A} (f), \ \ \overline g | \ ] \tag 3.3
$$
of
$$
\widetilde\omega_{\varepsilon} (f) \ \in \ \Omega_{m - n} (E (f);
\widetilde\varphi) = \bigoplus_{A \in \pi_0 (E (f))} \Omega_{m -
n} (A; \widetilde\varphi | A) \tag 3.4
$$
is nonzero.

The (nonnegative, integer) \ {\it Nielsen number}
$$
N_{\varepsilon} (f) \ := \ \# \{ A \in \pi_0 (E (f)) | \
\widetilde\omega_{\varepsilon, A} (f) \ne 0 \} \tag 3.5
$$
counts these essential pathcomponents of \ $E (f)$ \ (or,
equivalently, the essential Reidemeister classes). This is a
refinement of the concept of Nielsen numbers studied in \ \cite{K
5}, and we have
$$
N (f'_1 \scirc p_1, f'_2 \scirc p_2) \ \le \ N_+ (f) + N_- (f) \ .
$$

Clearly \ $N_{\varepsilon} (f)$ vanishes if and only if \
$\omega_{\varepsilon} (f)$ \ does. In other respects  the Nielsen
number is much cruder than the invariant \
$\widetilde\omega_{\varepsilon} (f)$ \ which, however, has the
drawback that it lies in a group which varies with \ $f$.

\proclaim{Proposition 3.6} \ If two link maps are link homotopic,
then they have the same Nielsen numbers.
\endproclaim

\demo{Proof} \ According to \ \cite{K 5} (see the discussion of
4.4) \ any homotopy \ $F$ \ from \ $f$ \ to another link map \
$\widehat f$ \ yields a homotopy equivalence \ $E (f) \sim E
(\widehat f)$ \ and an isomorphism of normal bordism groups which
maps the full coincidence invariant \ $\widetilde\omega (f'_1
\scirc p_1, f'_2 \scirc p_2)$ \ to \ $\widetilde\omega (\widehat
f'_1 \scirc p_1, \widehat f'_2 \scirc p_2)$. \ If \ $F$ \ is a
link homotopy this isomorphism preserves also the decomposition \
$\widetilde\omega = \widetilde\omega_+ + \widetilde\omega_-$.
\hfill $\blacksquare$
\enddemo

 The Nielsen number \ $N_{\varepsilon} (f), \ \varepsilon = + $ \
or \ $- \ $, \ is an example of a link homotopy invariant
extracted from \ $\widetilde\omega_{\varepsilon} (f)$ \ and taking
values in a set which is independent of \ $f$. Another such
example is \ $\omega_{\varepsilon} (f)$ \ (in the special case
where \ $N$ \ is stably parallelizable and hence \ $\varphi = -
TM$, \ cf.\ 2.10 and 2.8). The invariant \ $\alpha_w (f)$ \ of
Alexander Pilz (cf.\ \cite{P}, 3.9) can be interpreted as a third
such example: \ assume \ $f$ \ preserves base points (e.g.\ $f_i
(x_{i 0}) = (y_0, (- 1)^i), \ i = 1, 2)$ \ so that the
contributions \ $\widetilde\omega_{\varepsilon, A} (f)$ \ of the
various pathcomponents  \ $A \in \pi_0 (E (f))$ \ to \
$\widetilde\omega_{\varepsilon} (f)$ \ (cf.\ 3.3 and 3.4) can be
parameterized by the Reidemeister set \ $R$ \ (cf.\ 3.2) which
remains unchanged by base point preserving homotopies; then
$$
\alpha_w (f) \ = \ \{ [C_{+, A} (f)] \} \ \in \ \bigoplus_{A \in
R} \ \Omega_{m -n} \tag 3.7
$$
where \ $\Omega_*$ \ stands for framed or (un)oriented bordism
according as \ $M_1, M_2$ \ and \ $N$ \ are framed or
(un)oriented, resp. Thus \ $\alpha_w (f)$ neglects e.g.\ the map \
$\widetilde g | C_{+, A} (f)$ \ but still registers the
decomposition of the overcrossing locus \ $C_+ (f)$ \ into various
\lq\lq Reidemeister (or Nielsen) classes\rq\rq. \ If \ $m = n$, \
this enables A.\ Pilz to obtain full classification results (cf.\
\cite{P}, 5.4).
\medskip

Next we recall a result of Hatcher and Quinn (cf.\ \cite{HQ}, 3.1)
which allows  us sometimes to interpret \
$\widetilde\omega_{\varepsilon} (f)$ \ itself (without any loss of
information) as a link homotopy invariant.

\proclaim{Proposition 3.8} \ Assume that \ $M_1$ \ and \ $M_2$ \
are \ $(m - n + 1)$-connected.

Then the map \ $\incl$ \ $($cf.\ $3.1)$ -- together with a choice
of an orientation of \ $\varphi$ \ at \ $x_0$ -- induces an
isomorphism
$$
\incl_* \ : \ \Omega^{fr}_{m - n} (\Lambda (N, y_0)) \ @>{ \ \cong
\ }>> \ \Omega_{m -n} (E (f); \widetilde\varphi) \ .
$$

The framed bordism class
$$
\incl^{- 1}_* \ (\widetilde\omega_{\varepsilon} (f)) \ \in \
\Omega^{fr}_{m -n} (\Lambda (N, y_0)) \ , \ \ \ \varepsilon = +  \
\text{or} \ - ,
$$
is invariant under base point preserving link homotopies.
\endproclaim

\demo{Proof} \ Our choices (including those which are incorporated
in 3.1) induce isomorphisms
$$
\Omega^{fr}_{m -n} (\Lambda (N, y_0)) \ \cong \ \Omega^{fr}_{m -n}
(pr^{- 1} (x_0)) \ \cong \ \Omega_{m -n} (pr^{- 1} (x_0);
\widetilde\varphi |)
$$
and so does the fiber inclusion; this follows via a cell-by-cell
argument applied to (projected) maps into \ $M$.

Similarly if \ $F$ \ is a base point preserving link homotopy, the
(generic) coincidence manifold \ $C_{\varepsilon} (F) \ \subset \
M \times I$ \ can be retracted to \ $\{ x_0\} \times I$ \ and
hence yields the required bordism in the fiber \ $pr^{- 1} (x_0)
\sim \Lambda (N, y_0)$. \hfill $\blacksquare$
\enddemo

\remark{Remark {\rm 3.9}} \ If  \ $n > 0$, then every (free) link
homotopy class can be represented by a base point preserving link
map \ $f$ (i.e. $f_1 (x_{01}) \ne f_2 (x_{02})$ \ are fixed
preassigned points in \ $N \times \Bbb R$).

If two such link maps \ $f$ \ and \ $\widehat f$ \ are related by
a free link homotopy \ $F$, \ then we have (under the assumptions
of proposition 3.8 and provided \ $m_1, m_2 \le n - 2$) that
$$
\incl^{- 1}_* (\widetilde\omega_{\varepsilon} (\widehat f\ )) \ =
\ \sign (\tau_1) \cdot \sign (\tau_2) \cdot  c (\tau_1, \tau_2)_*
(\incl^{- 1}_* (\widetilde\omega_{\varepsilon} (f)))
$$
where \ $\tau_i = F'_i (x_{0i}, -)$ \ denotes the loop traced out
in \ $N$ \ during the homotopy, \ $i = 1,2$, \ and the
self-homotopy equivalence \ $c (\tau_1, \tau_2)$ \ of \ $\Lambda
(N, y_0)$ \ is defined by
$$
c (\tau_1, \tau_2)(\rho) \ = \ \sigma_1 \tau^{- 1}_1 \sigma^{-
1}_1 \rho \sigma_2 \tau_2 \sigma^{- 1}_2 \ \ \ , \ \ \ \ \rho \in
\Lambda (N, y_0) \ ;
$$
moreover, \ $\sign (\tau_i)$ \ equals \ $+ 1$ \ or \ $- 1$ \
according as \ $\tau_i$ \ preserves the orientation of \ $N$ \ or
not.
\endremark

\vskip7mm

\specialhead  4.\ \ Products of spheres
\endspecialhead

In this section we consider the special case
$$
N \ =  \ (S^1)^q \times S^2 \times \dots \times S^{r_\ell} \ = \
(S^1)^q \times N'
$$
where
$$
N' = \ \prod^\ell_{i = 1} S^{r_i} \ , \qquad \ \ \ 2 \le r_1 \le
r_2 \le \dots \le r_\ell, \ \ 1 \le \ell < \infty \
$$
denotes the subproduct formed by the simply connected factor
spheres. Pick base points \ $*_1 = 1 \in S^1; \ *_{r_i} \in
S^{r_i}, \ i = 1, \dots, \ell; \ y'_0 = (*_{r_1}, \dots,
*_{r_\ell}) \in N'$ \ and \ $y_0 = (1, \dots, 1, y'_0) \in N$. \
In view of proposition 3.8 the following sample calculation
determines frequently the group in which our \
$\widetilde\omega$-invariants lie.

\proclaim{Theorem 4.1} \ There are canonical isomorphisms of \
$\Omega^{fr}_*$-modules
$$
\Omega^{fr}_* (\Lambda (N, y_0)) \ \cong \ \bigoplus_{j \in \Bbb
Z^q} \ \Omega^{fr}_* (\Lambda (N', y'_0))
$$
and
$$
h \ = \ \oplus h_k \ : \ \Omega^{fr}_* (\Lambda (N', y'_0)) \ @>{
\ \cong \ }>> \ \bigoplus_{k \in \Bbb N^\ell} \Omega^{fr}_{* - d
(k)}
$$
where \ $d (k) := \sum^\ell_{i = 1} k_i (r_i -1)$ \ for \ $k =
(k_1, \dots, k_\ell) \in \Bbb N^\ell$ \ and \ $\Bbb N$ \ denotes
the set of natural numbers including \ $0$.
\endproclaim

Thus the Pontryagin-Thom isomorphism  \ $\Omega^{fr}_* \cong
\pi^S_*$ \ makes the framed bordism groups of the loop space of \
$N$ \ accessible to standard methods of stable homotopy theory.

\demo{Proof} \ Consider \ $\Bbb Z^q$ \ as a discrete topological
subspace of \ $\Bbb R^q$ \ and apply the exponential map \ $\Bbb
R^q \to (S^1)^q$ \ to the straight path in \ $\Bbb R^q$ \ joining
\ $0$ \ to any point in \ $\Bbb Z^q$. This yields the homotopy
equivalences \ $\Bbb Z^q \sim \Lambda ((S^1)^q)$ \ and
$$
\Lambda (N, y_0) \ \sim \ \Bbb Z^q \times \Lambda (N', y'_0)
$$
as well as the first isomorphism claimed above.
\enddemo

The second isomorphism generalizes a geometric construction which
was discussed in detail in  \cite{K 5} and which is closely
related to James-Hopf-invariants.

Given a framed bordism class \ $\nu$ \ of \ $\Lambda (N', y'_0)$,
consider the adjoint
$$
v = (v_1, \dots, v_\ell) : (V \times \Bbb R \ \cup \{ \infty\}, \
\infty) \ \longrightarrow \ (N' = S^{r_1} \times \dots \times
S^{r_\ell},  \ y'_0)
$$
of a generic representative. Then for \ $i = 1, \dots, \ell$ \ the
inverse image \ $W_i := v_i^{- 1} (*'_i)$ \ of a point \ $*'_i \ne
*_{r_i}$ \ in \ $S^{r_i}$ \ is a smoothly embedded framed
submanifold  of \ $V \times \Bbb R$. As in \cite{K 5}, \S\ 8 (see
also \cite{KS}) we may even assume that this embedding projects to
a generic framed \ $(r_i - 1)$-codimensional immersion \ $e'_i$ \
into $V$. Its \ $k_i$-fold selfintersection yields a \ $k_i (r_i -
1)$-codimensional immersion \ $e^{k_i}_i$ which is again framed
since the intersection branches are ordered by the \ $\Bbb
R$-component of the embedding \ $W_i \subset V \times \Bbb R$. Now
for \ $k = (k_1, \dots, k_\ell) \in \Bbb N^\ell$ \ define \ $h_k
(\nu)$ \ to be the framed bordism class of the transverse
intersection of the immersions \ $e^{k_1}_1, \dots,
e^{k_\ell}_\ell$ \ in \ $V$.

A straightforward generalization of the proofs in \cite{K 5}, \S\
8 (where the case \ $\ell = 1$ \ was discussed) shows that \ $h_k$
\ is well-defined and \ $h = \oplus h_k$ \ is bijective. \hfill
$\blacksquare$

\remark{Remark $4.2$} \ The theorem above leads to an
interpretation in terms of Laurent polynominals:
$$
\split \Omega^{fr}_* (\Lambda (N, y_0)) \ & \cong \bigoplus_{(j_1,
\dots, j_q) \in \Bbb Z^q} \Omega^{fr}_* (\Lambda(N', y'_0))
X^{j_1}_1 \dots X^{j_q}_q \\ &= \ \ \ \ \ \Omega^{fr}_*
(\Lambda(N', y'_0)) [X^{\pm 1}_1, \dots, X^{\pm 1}_q] .
\endsplit
$$
Expressed in this language the operation \ $c (\tau_1, \tau_2)_*$
\ discussed in remark 3.9 amounts just to multiplication with a
fixed monomial \ $X^{j_1}_1 \dots X^{j_q}_q$. Thus we conclude
(under the assumptions of proposition 3.8 and provided \ $m_1, m_2
\le n - 2$) that \ $\incl^{- 1}_* (\widetilde\omega_\varepsilon
(f))$, \ {\it considered up to multiplication with such
monomials}, \ is invariant under \ {\it base point free} \ link
homotopies.

(Note the analogy to Alexander polynomials).
\endremark

\vskip7mm

\specialhead  5.\ \ Spherical link maps
\endspecialhead

In this section we discuss the case where \ $M_i$ \ is the sphere
\ $S^{m_i}$, equipped with a base point \ $x_{0i}, \ i = 1, 2$.
Also choose base points \ $y_{01} \ne y_{02}$ \ and \ $y_0$ \ in \
$N$ \ and paths \ $\sigma_1$ \ and \ $\sigma_2$ \ in \ $N$ \
joining \ $y_0$ \ to \ $y_{01}$ \ and \ $y_{02}$, resp. Let \
$BLM_{m_1, m_2} (N \times \Bbb R)$ \ denote the set of base point
preserving link homotopy classes of base point preserving link
maps
$$
f  \ = \ f_1 \amalg f_2 \ : \ S^{m_1} \amalg S^{m_2} \
\longrightarrow \ N \times \Bbb R
$$
(i´.e.\ \ $f_i (x_{0i}) \ = \ (y_{0i}, 0))$, \ and define
$$
BLM^{(i)} \ := \ \{ [f] \in BLM_{m_1, m_2} (N \times \Bbb R) |  \
0 = [f_i] \in \pi_{m_i} (N \times \Bbb R)\} , \ \ \ \ i = 1, 2 .
\tag 5.1
$$

Now assume \ $1 \le m_1, m_2 \le n - 2$ \ throughout this section.
Then link maps as well as link homotopies can be deformed until
the $i$-th component avoids \ $\{y_{0i \pm 1}\} \times \Bbb R
\subset N \times \Bbb R, \ i = 1,2$. Thus we can add two link maps
\ $f$ \ and \ $\widehat f$ \ by \lq\lq stacking them on top of one
another\rq\rq \ w.r.\ to the height given by the \ $\Bbb
R$-component in \ $N \times \Bbb R$: shift the link map \
$\widehat f$ \ in the positive \ $\Bbb R$-direction until it is
entirely above \ $f$ \ and join base points along \ $\{ y_{0i}\}
\times \Bbb R$. The resulting addition makes \ $BLM_{m_1, m_2} (N
\times \Bbb R)$ \ into a semigroup with null element.

According to proposition 3.8 our choices determine isomorphisms
allowing us to identify our \ $\widetilde\omega$-invariants with
elements in a \ {\it fixed} \ group which does not vary with \ $f$
\ anymore (for a more explicit description see the proof below).
This greatly simplifies statements about link homotopy invariance,
additivity, value sets etc.

\proclaim{Proposition 5.2} \ Assume \ $1 \le m_1, m_2 \le n - 2$.

Then for \ $\varepsilon = +$ \ and \ $\varepsilon = -$ \ the \
$\widetilde\omega_\varepsilon$-obstruction determines a
well-defined map
$$
\widetilde\omega_\varepsilon \ : \ BLM_{m_1, m_2} (N \times \Bbb
R) \ \longrightarrow \ \Omega^{fr}_{m - n} (\Lambda (N, y_0)) \ .
$$
For all \ $[f], [{\widehat f}] \in BLM_{m_1, m_2} (N \times \Bbb
R)$ \ we have
$$
\widetilde\omega_+ ([f] + [\widehat f]) \ = \ \widetilde\omega_+
([f]) + \widetilde\omega_+ ([\widehat f]) + \widetilde\omega (f'_1
\circ p_1, \widehat f'_2 \circ p_2)
$$
and
$$
\widetilde\omega_- ([f] + [\widehat f]) \ = \ \widetilde\omega_-
([f]) + \widetilde\omega_- ([\widehat f]) + \widetilde\omega
(\widehat f'_1 \circ p_1, f'_2 \circ p_2)
$$
where the last summand denotes the coincidence invariant
$($compare $2.9$ and \cite{K 5}, $4.4$ and $2.4)$ of the indicated
projections to \ $N$.

In particular, \ $\widetilde\omega_\varepsilon | BLM^{(i)}$ \ is a
homomorphism of semigroups; its image is a group and invariant
under both the left and the right action of \ $\pi_1 (N)$ \  on \
$\Omega^{fr}_{m -n} (\Lambda (N, y_0)), \ i = 1, 2$.
\endproclaim

\demo{Proof} \ Let us first specify a representative \ $\widetilde
g_\varepsilon  :  C_\varepsilon (f)  \to  \Omega_{m - n} (\Lambda
(N, y_0))$ \ of \ $\widetilde\omega_\varepsilon (f)$. Since \ $m -
n + 1 < m_i$ \ the projection of \ $C_\varepsilon (f) \subset
S^{m_1} \times S^{m_2}$ \ to \ \ $S^{m_i}$ \ generically avoids
the (antipodal) point \ $- x_{0i}$ \ and hence allows a retraction
\ $\rho_i$ \ in \ $S^{m_i}$ \ to the base point \ $x_{0i}, \ i =
1, 2$. For any \ $x \in C_\varepsilon (f)$ \ the loop \
$\widetilde g_\varepsilon (x)$ \ can then be described as the
composite of paths
$$
y_0 \ @>{ \ \sigma_1 \ }>> \ y_{01} \ @>{ f'_1 (\rho_1 (x,
-))^{-1}}>> \  f'_1 (x) =  f'_2 (x) \ @>{f'_2 \rho_2 (x, -)}>> \
y_{02} \ @>{ \ \sigma^{- 1}_2 \ }>> \ y_0 \ .
$$

If we compose \ $f$ \ with base point preserving reflection \ $r$
\ of \ $S^{m_1}$ \ or \ $S^{m_2}$ \ we change the framing of \
$C_\varepsilon (f)$ \ and its location in \ $S^{m_1} \times
S^{m_2}$ \ but not the homotopy class of \ $\widetilde
g_\varepsilon$. Thus
$$
\widetilde\omega_\varepsilon (f_1 \scirc r \amalg f_2) \ = \
\widetilde\omega_\varepsilon (f_1 \amalg f_2 \scirc r) \ = \ -
\widetilde\omega_\varepsilon (f) \ . \tag 5.3
$$

If \ $\tau_i$ \ is a loop in \ $N$ \ starting and ending at \
$x_{0i}$ \ (and generically avoiding \ $f_{i \pm 1}$) \ and we use
it to modify \ $f_i$ \ via the standard operation, we conclude for
the resulting link map that
$$
\widetilde\omega_\varepsilon (\tau_{1 *} (f_1) \amalg \tau_{2 *}
(f_2)) \ = \ \pm (\sigma_1 \tau^{- 1}_1 \sigma^{- 1}_1)
\widetilde\omega_\varepsilon (f) (\sigma_2 \tau_2 \sigma^{- 1}_2)
\ . \tag 5.4
$$

Finally note that the overcrossing locus \ $C_\varepsilon (f +
\widehat f)$ \ consists of \ $C_\varepsilon (f) \cup C_\varepsilon
(\widehat f)$ \ and of the \ {\it full} \ coincidence locus of \
$(f'_1, \widehat f'_2)$. \ The proposition follows. \hfill
$\blacksquare$
\enddemo

Next consider the invariant
$$
\widetilde{\widetilde \omega}_\varepsilon \ : \ BLM_{m_1, m_2} (N
\times \Bbb R) \ \longrightarrow \ \Omega^{fr}_{m -n} (\Lambda (N,
y_0)) \times \bigoplus^2_{i = 1} \pi_{m_i} (N, y_{0i})
$$
which enriches \ $\widetilde\omega_\varepsilon (f)$ \ by the
homotopy classes \ $[f'_1]$ \ and \ $[f'_2]$ \ of the component
maps \ $f_1$ \ and \ $f_2$ \ of \ $f$. In view of 2.9 \ \
$\widetilde{\widetilde\omega}_+$ \ determines \
$\widetilde{\widetilde\omega}_-$.

\proclaim{Corollary 5.5} \ Both for \ $\varepsilon = +$ \ and \
$\varepsilon = -$ \ we have
$$
\widetilde{\widetilde\omega}_\varepsilon (BLM_{m_1, m_2} (N \times
\Bbb R)) \ = \ \widetilde\omega_\varepsilon (BLM^{(1)} \cap
BLM^{(2)}) \times \bigoplus^2_{i = 1} \pi_{m_i} (N, y_{0i})
$$
$($compare $5.1)$.
\endproclaim

\demo{Proof} \ Given \ $i = 1, 2$, \ let
$$
c_i \ \ : \ \ \pi_{m_i} (N, y_{0i}) \ \cong \ \pi_{m_i} (N -
\{y_{0i \pm 1}\}, y_{0i}) \ \longrightarrow \ BLM^{(i)}
$$
be defined by adding the constant map with value \ $y_{0 i \pm 1}
\in N = N \times \{0\} \subset N \times \Bbb R$. \ Then we have
the bijection
$$
(BLM^{(1)} \cap BLM^{(2)}) \ \times \ \bigoplus^2_{i = 1}
\pi_{m_i} (N, y_{0i}) \ \longrightarrow BLM_{m_1, m_2} (N \times
\Bbb R) , \tag 5.6
$$
$(f, [\widehat f_1], [\widehat f_2]) \ \longrightarrow \ c_2
([\widehat f_2]) + f + c_1 ([\widehat f_1])$. \ Here the order of
summation is chosen in such a way that according to 5.2 the sum
has the same \ $\widetilde\omega_+$-value as \ $f$. \ The
corollary follows. \hfill $\blacksquare$
\enddemo

\demo{Proof of theorem 1.1 of the introduction} \ Clearly the
claim holds if \ $m_1 = 0$.

Thus assume that \ $m_1 \ge 1$. \ Then also \ $m_1, m_2 \leq n-2$\
. In view of the decomposition 5.6 we have to prove our claim only
for \ $f, \widehat f \in BLM^{(1)} \cap BLM^{(2)}$ \ with \
$\widetilde\omega_+ (f) = \widetilde\omega_+ (\widehat f)$. \ Put
\ $(- \widehat f) := \widehat f_1 \scirc r \amalg \widehat f_2$ \
(cf.\ 5.3). Then according to 5.2 we have
$$
\widetilde\omega_+ (f + (- \widehat f)) \ = \ \widetilde\omega_+
(f) - \widetilde\omega_+ (\widehat f) \ = \ 0 \ .
$$
Thus it follows from theorem 2.13 that \ $f + (- \widehat f)$ \
(and similarly \ $(- \widehat f) + \widehat f\ )$ \ is link
nulhomotopic. Therefore
$$
f \ \ \sim \ f \ + ((- \widehat f) + \widehat f) \ \sim \ (f + (-
\widehat f)) \ + \ \widehat f \ \ \sim \ \ \widehat f \ .
$$
\vskip-5mm \hfill $\blacksquare$
\enddemo

\vskip4mm Finally  we discuss a simple construction which produces
many link maps with interesting \
$\widetilde\omega_\pm$-invariants.

Consider the case where \ $N = N_1 \times N_2$ \ is the product of
two manifolds of dimensions \ $n_1, n_2 \ge 1$. Then the base
points take the form \ $y_{0i} = (y^1_{0i}, y^2_{0i}), \ i = 1,2$,
and we may assume that \ $y^1_{01} \ne y^1_{02}$ \ and \ $y^2_{01}
\ne y^2_{02}$.

Also consider the link map
$$
e \ = \ e_1 \amalg e_2 \ : \ S^{n_1} \amalg N_2 \ \longrightarrow
\ N \times \Bbb R \tag 5.7
$$
formed by the embedding
$$
e_2 \ : N_2 \ \cong \ \{y^1_{02}\}  \times N_2 \times \{0\} \
\subset \ N \times \Bbb R
$$
and by a meridian (of radius $\rho$, say) which is based at \
$(y_{01},0)$ \ and which lies in \ $N_1 \times \{y^2_{01}\} \times
\Bbb R$ \ (so that \ $(y^1_{02}, y^2_{01})$ \ is the only
over/undercrossing point). Note that \ $e_1$ \ can be contracted
within a normal ball of \ $e_2$ \ with center at \ $(y^1_{02},
y^2_{01}, 0)$.

Thus composition with \ $e$ \ yields a map
$$
\pi_{m_1} (S^{n_1}) \ \times \ \pi_{m_2} (N_2) \ \longrightarrow \
BLM^{(1)} \tag 5.8
$$
(compare 5.1). Moreover the inclusion \ $e'_2$ \ induces a
monomorphism
$$
e'_{2 *} \ : \ \Omega^{fr}_* (\Lambda (N_2, y^2_{02})) \
\longrightarrow \ \Omega^{fr}_* (\Lambda (N, y_0)) \ .
$$

\proclaim{Proposition 5.9} \ Assume \ $1 \le m_1, m_2 \le n_1 +
n_2 - 2$ \ . Then we have for all \ $[f_1] \in \pi_{m_1}
(S^{n_1}), \ [f_2] \in \pi_{m_2} (N_2)$ \ and both for \
$\varepsilon = +$ \ and \ $\varepsilon = -$ \ that
$$
\widetilde\omega_\varepsilon (e_1 \scirc f_1 \ \amalg \ e_2 \scirc
f_2) \ \ = \ \ \pm \deg (f_1) \cdot e'_{2*} (\widetilde\deg (f_2))
$$
where the degrees
$$
\deg (f_1) \ := \ \omega (f_1, *) \ \in \ \Omega^{fr}_{m_1 - n_1}
$$
and
$$
\widetilde\deg (f_2) \ := \ \widetilde\omega (f_2, y^2_{01}) \ \in
\ \Omega^{fr}_{m_2 - n_2} (\Lambda (N_2, y^2_{02}))
$$
measure coincidences with the indicated constant maps $($as in
\cite{K 5}, {\rm 1.11, 1.12,} and {\rm 7.5)}.
\endproclaim

\demo{Proof} \ Define \ $L_1 := f_1^{-1}(\{(y^1_{02}, -\varepsilon
\rho)\})$ \ and $L_2 := f_2^{-1} (\{y^2_{01}\})$.\ Then the over/
undercrossing locus is
$$
C_\varepsilon \ = L_1 \times L_2 \ \subset \ S^{m_1} \times
S^{m_2} .
$$
The corresponding loops are described as in the proof of
proposition 5.2:  \ apply \ $f_i$ \ to the paths resulting from a
contraction of \ $L_i , \ i = 1, 2$. Since \ $f_1$ \ is
nulhomotopic, the framed manifold \ $L_1$ \ is equipped with a
trivial map and contributes to \ $\widetilde\omega_\varepsilon$ \
only via the \ $\Omega^{fr}_*$-module structure on \
$\Omega^{fr}_* (\Lambda (N, y_0))$. \hfill $\blacksquare$
\enddemo

Often it is convenient to identify the framed bordism ring \
$\Omega^{fr}_*$ \ with the stable homotopy ring \ $\pi^S_*$ \ of
spheres via the Pontryagin-Thom isomorphism. Then \ $\deg (f_1)$ \
corresponds to the stable Freudenthal suspension \ $\pm E^\infty
([f_1])$.

\proclaim{Corollary 5.10} \ Both for \ $\varepsilon = +$ \ and \
$\varepsilon = -$ \ the value set \ $\widetilde\omega_\varepsilon
(BLM^{(1)})$ \ $($cf.\ $5.1$ and $5.2)$ contains at least the
subgroup generated by the set
$$
E^\infty (\pi_{m_1} (S^{n_1})) \ \cdot \ e'_{2*}(\widetilde\deg
(\pi_{m_2} (N_2))) \ \ \subset \ \ \Omega^{fr}_{m - n} (\Lambda
(N, y_0))
$$
and by the left and right group action of \ $\pi_1 (N, y_0)$.
\endproclaim

This follows from propositions 5.2 and 5.9.

\example{Example {\rm 5.11} : Spherical link maps into products of
spheres} \ Let
$$
N \ = \ (S^1)^q \ \times \ S^2 \ \times \dots \times  S^{r_\ell}
$$
be as in section 4. Given \ $1 \le i \le \ell$, \ let \ $N_2$ \
and \ $N_1$ \ be the $i$-th higher dimensional sphere \ $S^{r_i}$
\, (i.e. \ $r_i \ge 2$) \, and the product of the remaining
factors, resp. Also recall from theorem 1.14 in \cite{K 5} that \
$\widetilde\deg (\pi_{m_2} (S^{r_i}))$ \ corresponds to the image
of the total stabilized James-Hopf-invariant homomorphism
$$
\Gamma_i = \bigoplus_{k \ge 0} E^\infty \scirc \gamma_{k+1} \ : \
\pi_{m_2} (S^{r_i}) \ \longrightarrow \ \bigoplus_{k  \ge 0}
\pi^S_{m_2 - r_i - k (r_i - 1)} \ \ .
$$
Thus according to corollary 5.10 at least all the elements of the
subgroup
$$
\bigoplus_{j \in \Bbb Z^q} \ E^\infty (\pi_{m_1} (S^{n - r_i})) \
\cdot \ \Gamma_i (\pi_{m_2} (S^{r_i})) \tag 5.12
$$
of
$$
\Omega^{fr}_{m - n} (\Lambda (N, y_0)) \ \cong \ \bigoplus_{j \in
\Bbb Z^q} \ \bigoplus_{k \in \Bbb N^\ell} \ \pi^S_{m - n - d (k)}
$$
(cf.\ theorem 4.1) occur as values of our invariants \
$\widetilde\omega_+$ \ and \ $\widetilde\omega_-$ \ for suitable
link maps, \ $i = 1, \dots, \ell.$
\endexample

\example{Special case 5.11a: $\bold N = \bold S^2 \bold\times
\bold S^3, \ \bold m_1 = \bold m_2 = 3$} \ \ Here \
$\widetilde\omega_+$ \ and \ $\widetilde\omega_-$ \ take values in
$$
\Omega^{fr}_1 (\Lambda N) \ \cong \ \Omega^{fr}_1 (\Lambda S^2) \
@>{ \ \cong \ }>{ \ h \ }> \ \Omega^{fr}_1 \oplus \Omega^{fr}_0 =
\Bbb Z_2 \oplus \Bbb Z
$$
(cf.\ 4.1). If we pick \ $N_2 = S^2$ \ and \ $N_1 = S^3$ \ we see
that
$$
(1, 1) = (E^\infty, \ \text{Hopf invariant) (Hopf map)} = \Gamma_1
\ \text{(Hopf map)}
$$
corresponds to a value of \ $\widetilde\omega_{\pm}$; if we pick \
$N_2 = S^3, N_1 = S^2$ \ we obtain the same for \ $(1 = E^\infty
\text{(Hopf map)}, 0)$. Thus the invariant \
$\widetilde\omega_\pm$ (even when restricted to \ $BLM^{(i)}, \ i
= 1$ or $2$, cf.\ 5.1) assumes all possible values in its target
group.

Both the invariant \ $\omega_+$ \ (cf.\ 2.10) and the Pilz
invariant \ $\alpha_w$ \ (cf. 3.7) capture only the \ $\Bbb
Z_2$-component of \ $\widetilde\omega_+$, i.e.\ the framed bordism
class of the overcrossing locus \ $C_+$ \ (without any nontrivial
map into a target space or any nontrivial decomposition into
Nielsen classes).
\endexample

\example{Special case 5.11b : $\bold N = (\bold S^1)^3 \bold\times
\bold S^2, \  \bold m_1 = \bold m_2 = 3$} \ \ By theorem 4.1 the
target group of \ $\widetilde\omega_\pm$ \ is as stated in the
example 1.3 of the introduction. Moreover, every element \
$(1,1)_j, \ j \in \Bbb Z^3$, \ lies in \ $\widetilde\omega_\pm
(BLM^{(i)}), \ i = 1, 2$ (cf.\ 5.1, 5.2, and 5.11a). So does also
\ $(1,0)_j, \ j \in \Bbb Z^3$. \ Indeed, the Fenn-Rolfsen link map
\ $f_{FR} : S^2 \amalg S^2 \to \Bbb R^4$ \ (cf.\ \cite{FR}) can be
suspended w.r.\ to each component to yield a link map \ $\Sigma_1
(\Sigma_2 (f_{FR})) : S^3 \amalg S^3 \to \Bbb R^6$ \ with
nontrivial \ $\alpha$-invariant (cf. e.g.\ \cite{K 2}, 2.1, 2.10,
and 4.8); included in a small ball in \ $N \times \Bbb R$ \ and
connected to the base points in a suitable way this yields a link
map in \ $BLM^{(1)} \cap BLM^{(2)}$ \ with the desired \
$\widetilde\omega_\pm$-value. Thus
$$
\widetilde\omega_\varepsilon (BLM^{(i)}) \ = \ \bigoplus_{j \in
\Bbb Z^3} \ (\Bbb Z_2 \oplus \Bbb Z)
$$
for every combination of \ $i = 1, 2$ and $\varepsilon = +$ or
$-$.

In this special case the Pilz invariant \ $\alpha_w \in \oplus
\Bbb Z_2$ \ is obtained by projecting to the \ $\Bbb Z_2$-factors.
It is much stronger here than the invariant \ $\omega_+ \in \Bbb
Z_2 \cong \Omega^{fr}_1 \cong \Omega^{fr}_1 (S^3 \times S^3)$ \
which measures just the framed bordism class of the overcrossing
locus without retaining its Nielsen decomposition.
\medskip

On the other hand in some situations Nielsen decompositions are
irrelevant and \ $\omega_+$ \ captures more information than the
Pilz invariant \ $\alpha_w$.
\endexample

\example{Example 5.13} \ Given an integer \ $s \ge 1$, \ equip \
$S^s = \partial(B^{s +1})$ \ with the boundary framing inherited
from the \ $(s + 1)$-ball and consider the link map
$$
f \ = \ f_1 \amalg f_2 \ : \ S^3 \times S^s \amalg S^3 \
\longrightarrow \ S^2 \times S^3 \times \Bbb R
$$
defined by \ $f_1 (x_1, x_2) := (\text{(Hopf map)} (x_1), *_{S^3},
0)$ \ for \ $(x_1, x_2) \in S^3 \times S^s$ \ and \ $f_2 (x) :=
(*_{S^2}, x, 1)$ \ for $x \in S^3$. \ Then \ $\alpha_w (f) = 0$ \
while \ $\omega_+ (f)$ \ is nontrivial. Indeed, the coincidence
locus \ $C_+ = S^1 \times S^s$ \ is framed nulbordant, but
$$
\omega_+ (f) \ = \ [C_+, g_+, \overline g_+] \ \in \
\Omega^{fr}_{s + 1} (S^3 \times S^s \times S^3) .
$$
projects to
$$
(0, 1) \ \in \ \Omega^{fr}_{s +1} (S^s) \ \cong \ \Omega^{fr}_{s +
1} \oplus \Omega^{fr}_1
$$
(compare 2.15).
\endexample

\vskip1cm

\head{References}
\endhead

\enddocument